\documentclass[10pt]{article}

\usepackage{amsmath}
\usepackage{amsfonts}
\usepackage{amssymb}
\usepackage[english]{babel}
\usepackage{amsthm}
\def\pf{\par\noindent {\bf Proof}~\par\noindent}

\newcommand{\mR}{\mathbb{R}}
\newcommand{\mC}{\mathbb{C}}
\newcommand{\mN}{\mathbb{N}}

\newcommand{\mE}{\mathbb{E}}
\newcommand{\mS}{\mathbb{S}}
\newcommand{\mH}{\mathbb{H}}
\newcommand{\mI}{\mathbb{I}}
\newcommand{\mJ}{\mathbb{J}}
\newcommand{\mK}{\mathbb{K}}

\newcommand{\mV}{\mathbb{V}}
\newcommand{\mW}{\mathbb{W}}

\newcommand{\mcP}{\mathcal{P}}

\newcommand{\mcH}{\mathcal{H}}

\newcommand{\mcE}{\mathcal{E}}

\newcommand{\gf}{\mathfrak{f}}
\newcommand{\gfd}{\mathfrak{f}^{\dagger}}
\newcommand{\gsl}{\mathfrak{sl}}
\newcommand{\gsp}{\mathfrak{sp}}
\newcommand{\gog}{\mathfrak{g}}

\newcommand{\gso}{\mathfrak{so}}
\newcommand{\ggl}{\mathfrak{gl}}

\newcommand{\gu}{\mathfrak{u}}

\newcommand{\ol}{\overline}
\newcommand{\olz}{\ol{z}}
\newcommand{\onehalf}{\frac{1}{2}}
\newcommand{\HS}{\mcH^S}
\newcommand{\HSD}{\mcH^{S^{\dagger}}}
\newcommand{\SE}{\mcE}
\newcommand{\SED}{\mcE^{\dagger}}
\newcommand{\domcodom}{(\mR^{4p};\mC)}
\newcommand{\eop}{\hfill$\square$}

\newcommand{\uz}{\underline{z}}
\newcommand{\uzj}{\underline{z}^J}
\newcommand{\uzd}{\underline{z}^{\dagger}}
\newcommand{\uzdj}{\underline{z}^{\dagger J}}

\newcommand{\upz}{\partial_{\uz}}
\newcommand{\upzj}{\partial_{\uz}^J}

\newcommand{\upzd}{\partial_{\uz}^{\dagger}}
\newcommand{\upzdj}{\partial_{\uz}^{\dagger J}}

\newcommand{\nz}{\vert\underline{z}\vert}

\newcommand{\p}{\partial}

\newtheorem{theorem}{Theorem}
\newtheorem{lemma}{Lemma}
\newtheorem{proposition}{Proposition}
\newtheorem{remark}{Remark}
\newtheorem{corollary}{Corollary}

\newtheorem{definition}{Definition}

\begin{document}

\title{Fundaments of Quaternionic Clifford Analysis III: \\
 Fischer Decomposition in Symplectic Harmonic Analysis}
\author{F.\ Brackx$^\ast$, H.\ De Schepper$^\ast$, D.\ Eelbode$^{\ast\ast}$, R.\ L\'{a}vi\v{c}ka$^\ddagger$ \& V.\ Sou\v{c}ek$^\ddagger$}

\date{\small{$^\ast$ Clifford Research Group, Faculty of Engineering and Architecture, Ghent University\\
Building S22, Galglaan 2, B-9000 Gent, Belgium\\
$\ast\ast$ University of Antwerp, Middelheimlaan 2, Antwerpen, Belgium\\
$^\ddagger$ Charles University in Prague, Faculty of Mathematics and Physics, Mathematical Institute\\
Sokolovsk\'a 83, 186 75 Praha, Czech Republic}}

\maketitle


\begin{abstract}
In the framework of quaternionic Clifford analysis in Euclidean space $\mR^{4p}$, which constitutes a refinement of Euclidean and Hermitian Clifford analysis, the Fischer decomposition of the space of complex valued polynomials is obtained in terms of spaces of so--called (adjoint) symplectic spherical harmonics, which are irreducible modules for the symplectic group Sp$(p)$. Its Howe dual partner is determined to be $\gsl(2,\mC) \oplus \gsl(2,\mC) =  \gso(4,\mC)$.
\end{abstract}


\section{Introduction}

In 1917 Ernst Fischer proved (see \cite{fischer}) that, given a homogeneous polynomial $q(X)$, $X \in \mR^m$, every homogeneous polynomial $P_k(X)$ of degree $k$ can be uniquely decomposed as 
$$
P_k(X) = Q_k(X) + q(X)R(X)
$$ 
where $Q_k(X)$ is a homogeneous polynomial of degree $k$ satisfying the partial differential equation
$$
q(D)Q_k = 0
$$
$D$ being the differential operator corresponding to $X$ through Fourier identification ($X_j \leftrightarrow \p_{x_j}, j=1, \ldots,m$) and $R(X)$ is a homogeneous polynomial of suitable degree. If, in particular, $q(X) = |X|^2 = \sum_{j=1}^m X_j^2 = r^2$, then $q(D) = \sum_{j=1}^m  \p_{X_j}^2 = \Delta_m$, the Laplace operator in $\mR^m$, and $Q_k$ is  harmonic, leading to the well-known decomposition 
\begin{equation}
\label{harm-decomp}
\mcP(\mR^m;\mC) = \bigoplus_{k=0}^{\infty} \bigoplus_{p=0}^{\infty} r^{2p} \ \mcH_k(\mR^m;\mC)
\end{equation}
of the space $\mcP(\mR^m;\mC)$ of complex valued polynomials, in terms of the spaces $\mcH_k(\mR^m;\mC)$ of complex valued harmonic homogeneous polynomials of degree $k$.
This space $\mcP(\mR^m;\mC)$ is a module over the special orthogonal group SO($m$), its action being the regular representation
\begin{equation}
\label{actionSO}
[ g \cdot P ](X)= P(g^{-1} \cdot X), \quad g \in \mbox{SO}(m), \quad P \in \mcP(\mR^m;\mC), \quad X \in \mR^m
\end{equation}
Each of the constituents of the decomposition (\ref{harm-decomp})
$$
r^{2p} \ \mcH_k, \quad p \in \mN_0:=\mN \cup \{0\}, \; \quad k \in \mN_0
$$
is a subspace of $\mcP(\mR^m;\mC)$ which is invariant under the SO($m$)--action and all SO($m$)-modules $\mcH_k(\mR^m;\mC)$ are irreducible and mutually inequivalent. In particular, the space  $\mcP_k(\mR^m;\mC)$ of homogeneous polynomials of degree $k$, decomposes into SO($m$)-irreducibles 
as
\begin{equation}
\label{homog-harm-decomp}
\mcP_k(\mR^m;\mC) = \bigoplus_{p=0}^{\lfloor  \frac{k}{2} \rfloor} r^{2p} \ \mcH_{k-2p}(\mR^m;\mC)
\end{equation}

\noindent
The Fischer decomposition (\ref{harm-decomp}) may be rewritten in the triangular diagram
\begin{equation}
\label{harm-decomp-triangular}
\begin{array}{ccccccc}
\mcH_0 && r^2 \ \mcH_0 && r^4 \ \mcH_0 && \cdots\\
& \mcH_1 && r^2 \ \mcH_1 && \cdots \\
& & \mcH_2 && r^2 \ \mcH_2 && \cdots \\
& & & \mcH_3 && \cdots \\
& & & & \mcH_4 && \cdots
\end{array}
\end{equation}
the vertical columns then reflecting the decomposition (\ref{homog-harm-decomp}) of the spaces $\mcP_k(\mR^m;\mC)$, $k=0,1,2,\ldots$\\[-2mm]

\noindent
It is  clear that in the Fischer decompositions (\ref{harm-decomp}) and (\ref{homog-harm-decomp}) the operators $X := \frac{1}{2} \ r^2$ and $Y := - \frac{1}{2} \ \Delta_m$ play a key role. Note that they correspond to each other under natural or Fourier duality, also known as Fischer duality. They both commute with the action (\ref{actionSO}) of SO($m$) on functions and on polynomials in particular, and their mutual commutator is
$$
\left[ X , Y \right] = \left[ \frac{1}{2} \ r^2 , - \frac{1}{2} \ \Delta_m \right] = \mE + \frac{m}{2}
$$
where $\mE = r \p_r = \sum_{j=1}^m X_j \p_{X_j}$ is the Euler operator in $\mR^m$. We then put 
$$
H := \mE + \frac{m}{2}
$$
and find that $[ H , X ] = 2X$ and $[ H , Y ] = -2Y$. This means that $\{H, X, Y \}$ generates a three--dimensional Lie algebra isomorphic with the Lie algebra $\gsl(2,\mC)$. The action of $\gsl(2,\mC)$ on the decompositions (\ref{harm-decomp}) and (\ref{homog-harm-decomp}) is:
\begin{eqnarray}
\label{action}
X : r^{2p} \ \mcH_k   & \longrightarrow &  r^{2p+2} \ \mcH_k \nonumber \\
Y : r^{2p} \ \mcH_k   & \longrightarrow &  r^{2p-2} \ \mcH_k\\
H : r^{2p} \ \mcH_k   & \longrightarrow &  r^{2p} \ \mcH_k \nonumber
\end{eqnarray}

Taking the dimension $m$ to be even: $m=2n$, the standard complex structure $\mI_{2n}$ on $\mR^{2n}$ is introduced as follows. Let $E_n$ denote the identity matrix in $M_n(\mC)$, the space of square $n \times n$ matrices with complex entries. Let
$$
\varphi_n : M_n(\mC) \longrightarrow M_{2n}(\mR)
$$
stand for the injective homomorphism embedding $M_n(\mC)$ into the space $M_{2n}(\mR)$ of square $2n \times 2n$ real matrices. This embedding may be realized by substituting for each complex entry $a + b i$, the $2 \times 2$ real matrix 
$\bigl(\begin{smallmatrix}
\phantom{-}a&b\\ -b&a
\end{smallmatrix} \bigr)$.
In $\mC^n$ multiplication by the imaginary unit $i$ is the $\mC$--linear transformation associated to the matrix $i E_n$. The standard complex structure $\mI_{2n}$ then is the complex linear real matrix 
$$
\mI_{2n} =  \varphi_n(i E_n) = {\rm diag}\bigl(\begin{smallmatrix}
\phantom{-}0&1\\ -1&0
\end{smallmatrix} \bigr)
$$
As expected, there holds $\mI_{2n}^2 = - E_{2n}$, $E_{2n}$ being the identity matrix in $M_{2n}(\mR)$. Moreover $\mI_{2n}$ belongs to SO$(2n)$, and a matrix $B \in M_{2n}(\mR)$ is complex linear, i.e. belongs to $\varphi_n(M_n(\mC))$, if and only if $B$ commutes with the complex structure $\mI_{2n}$ on $\mR^{2n}$. We then have the following result (see also \cite{FQA}).

\begin{proposition}
The {\rm SO}$(2n)$--matrices commuting with the complex structure $\mI_{2n}$ on $\mR^{2n}$ form a subgroup of  {\rm SO}$(2n)$, denoted by  {\rm SO}$_{\mI}(2n)$, which is isomorphic with the unitary group {\rm U}$(n)$.
\end{proposition}

\noindent
The introduction of the complex structure $\mI_{2n}$ allows for considering the space $\mcP(\mR^{2n};\mC)$ of complex valued polynomials defined on Euclidean space of {\em even} dimension, as an SO$_{\mI}(2n) \cong$ U$(n)$--module, the action of SO$_{\mI}(2n)$ being
$$
[ u \cdot P](X)= P(u^{-1} \cdot X), \quad u \in \mbox{SO}_{\mI}(2n), \quad P \in \mcP(\mR^{2n};\mC), \quad X \in \mR^{2n}
$$
Since each complex valued polynomial in the real variables $\left(X_1, \ldots, X_{2n}\right) = \left( x_1, \ldots, x_n, y_1, \ldots, y_n  \right)$ may be written as a polynomial in the complex variables $\left( z_1, \ldots, z_n, \olz_1, \ldots, \olz_n  \right)$, with $z_j=x_j + i \, y_j, \olz_j=x_j - i \, y_j,  j=1,\ldots,n$, i.e.
$$
P(X) = P(x_1, \ldots, x_n, y_1, \ldots, y_n) = \widetilde{P}(z_1, \ldots, z_n, \olz_1, \ldots, \olz_n)
$$
we have to determine the polynomials $\widetilde{P}$ which are invariant under the action of SO$_{\mI}(2n) \cong$ U$(n)$.
As is well--known, the space of U($n$)--invariant polynomials in $\mcP(\mR^{2n};\mC)$ is the space with basis
$$
\left( 1, r^2, r^4, \ldots, r^{2p}, \ldots   \right)
$$
where $r^2$ can be written as:
$$
r^2 = |X|^2 = \sum_{j=1}^{2n}  \ X_j^2 =  \sum_{j=1}^{n} \ x_j^2 + y_j^2 =  \sum_{j=1}^{n} \ z_j \olz_j =  \sum_{j=1}^{n} \ |z_j|^2
$$
The differential operator corresponding, under Fourier duality, to the generator $r^2$ is the Laplace operator 
$$
\Delta_{2n} = \sum_{j=1}^{n} \ \p^2_{x_j x_j} + \p^2_{y_j y_j} = 4 \ \sum_{j=1}^{n} \ \p_{z_j}\p_{\olz_j}
$$
whence we are led to consider the space of harmonic polynomials in $(z_1, \ldots, z_n, \olz_1, \ldots, \olz_n)$. Its subspace $\mcH_k(\mR^{2n};\mC)$ of complex valued $k$--homogeneous harmonic polynomials may be decomposed as
$$
\mcH_k(\mR^{2n};\mC) = \bigoplus_{a+b=k} \ \mcH_{a,b}(\mR^{2n};\mC)
$$
where $\mcH_{a,b}(\mR^{2n};\mC)$ is the space of the complex valued harmonic polynomials which are $a$--homogeneous in the variables $z_j$ and at the same time $b$--homogeneous in the variables $\olz_j$, i.e.
$$
H_{a,b}(\lambda z_1, \ldots, \lambda z_n, \mu \olz_1, \ldots, \mu \olz_n) = \lambda^a \ \mu^b \ H_{a,b}(z_1, \ldots, z_n, \olz_1, \ldots, \olz_n)
$$
This leads to the Fischer decomposition
\begin{equation}
\label{harmu}
\mcP(\mR^{2n};\mC) =
\bigoplus_{k=0}^{\infty} \bigoplus_{p=0}^{\infty}   \bigoplus_{a=0}^{k} \ r^{2p} \ \mcH_{a,k-a}(\mR^{2n};\mC)
\end{equation}
where the constituents
$$
r^{2p} \ \mcH_{a,k-a}, \quad p \in \mN_0, \quad k \in \mN_0, \quad a=0,\ldots,k
$$
are irreducible invariant subspaces under the action of U($n$). In particular, the space $\mcP_k(\mR^{2n};\mC)$ of $k$-homogeneous polynomials decomposes as
\begin{equation}
\label{homog-harmu}
\mcP_k(\mR^{2n};\mC) =
\bigoplus_{p=0}^{\lfloor  \frac{k}{2}\rfloor }  \bigoplus_{a=0}^{k-2p} \ r^{2p} \ \mcH_{a,k-2p-a}(\mR^{2n};\mC)
\end{equation}

\noindent
The corresponding diagram, similar to (\ref{harm-decomp-triangular}), looks like
\begin{equation}
\label{harmu-triangular}
\hspace*{-10mm}
\begin{array}{cccccccc}
\mcH_{0,0} && r^2 \ \mcH_{0,0} && r^4 \ \mcH_{0,0} && \cdots\\ \\
& \mcH_{1,0} && r^2 \ \mcH_{1,0} && r^4 \ \mcH_{1,0} && \cdots \\
& \mcH_{0,1} && r^2 \ \mcH_{0,1} && r^4 \ \mcH_{0,1} && \cdots \\ \\
& & \mcH_{2,0} && r^2 \ \mcH_{2,0} && \cdots \\
& & \mcH_{1,1} && r^2 \ \mcH_{1,1} && \cdots \\
& & \mcH_{0,2} && r^2 \ \mcH_{0,2} && \cdots \\ \\
& & & \mcH_{3,0}  && r^2 \mcH_{3,0} && \cdots \\
& & & \mcH_{2,1}  && r^2 \mcH_{2,1} && \cdots \\
& & & \mcH_{1,2}  && r^2 \mcH_{1,2} && \cdots \\
& & & \mcH_{0,3}  && r^2 \mcH_{0,3} && \cdots \\ \\
& & & & \mcH_{4,0} && \cdots\\
& & & & \mcH_{3,1} && \cdots\\
& & & & \mcH_{2,2} && \cdots\\
& & & & \mcH_{1,3} && \cdots\\
& & & & \mcH_{0,4} && \cdots\\�\\
& & & & &  \mcH_{5,0} && \cdots\\
& & & & &  \mcH_{4,1} && \cdots\\
& & & & &  \mcH_{3,2} && \cdots\\
& & & & &  \mcH_{2,3} && \cdots\\
& & & & &  \mcH_{1,4} && \cdots\\
& & & & &  \mcH_{0,5} && \cdots
\end{array}
\end{equation}

\noindent
The smallest Lie algebra of complex polynomial differential operators generated by the polynomial $r^2$ and its dual operator $\Delta_{2n}$ again is $\gsl(2,\mC)$, since
$$
\left[ X , Y \right] = \left[ \frac{1}{2} \ r^2 , - \frac{1}{2} \ \Delta_{2n} \right] = \mE + n = H
$$
However, there is an additional natural invariant differential operator coming into play. Indeed, the Euler operator
$\mE$ decomposes as 
$$\mE = \mE_z + \mE_{z}^{\dagger}$$ 
with
$$
\mE_z = \sum_{j=1}^{n}z_j \p_{z_j} \quad \mbox{and} \quad \mE_{z}^{\dagger} = \sum_{j=1}^{n}\olz_j\p_{\olz_j}
$$
Both these Euler operators in the complex variables are U$(n)$--invariant, and so is their difference, up to a chosen constant,
$$\mE_{z}^{\dagger} - \mE_{z} + n$$
which commutes with $X = \onehalf r^2, Y = - \onehalf \Delta_{2n}$ and $H =  \mE_z + \mE_{z}^{\dagger} + n$, since
$$
[r^2,\mE_z] =  - r^2, \quad [r^2,\mE_{z}^{\dagger}] = - r^2 
$$
and
$$
[\Delta_{2n},\mE_z] = \Delta_{2n}, \quad  [\Delta_{2n},\mE_{z}^{\dagger}] = \Delta_{2n}
$$
\noindent
In this way we end up with a reductive Lie algebra which is the direct sum of the three--dimensional Lie algebra generated by $\{H,X,Y\}$, isomorphic with $\gsl(2,\mC)$, and the one--dimensional abelian Lie algebra $\mC$ generated by $\{ \mE_{z}^{\dagger} - \mE_z + n \}$. This is nothing else but the four dimensional general linear algebra $\ggl(2,\mC)$ with action

\begin{eqnarray}
\label{uaction}
X : r^{2p} \ \mcH_{a,b}   & \longrightarrow &  r^{2p+2} \ \mcH_{a,b}\nonumber \\
Y : r^{2p} \ \mcH_{a,b}   & \longrightarrow &  r^{2p-2} \ \mcH_{a,b}\nonumber \\
H : r^{2p} \ \mcH_{a,b}   & \longrightarrow &  r^{2p} \ \mcH_{a,b}\\
\mE_{z}^{\dagger} - \mE_{z} + n: r^{2p} \ \mcH_{a,b}   & \rightarrow &  r^{2p} \ \mcH_{a,b}\nonumber
\end{eqnarray}

since also
\begin{eqnarray*}
\mE_{z}: r^{2p} \ \mcH_{a,b}   & \longrightarrow &  r^{2p} \ \mcH_{a,b}\\
\mE_{z}^{\dagger}: r^{2p} \ \mcH_{a,b}   & \longrightarrow &  r^{2p} \ \mcH_{a,b}
\end{eqnarray*}

When comparing the Fischer decompositions (\ref{harm-decomp}) and (\ref{harmu}), it becomes clear that refining the symmetry group from SO$(2n)$ to its subgroup SO$_{\mI}(2n) \cong$ U$(n)$, results into the splitting of the space $\mcH_k(\mR^{2n};\mC)$ of homogeneous harmonic polynomials, now considered as functions in the complex variables $\left( z_1, \ldots, z_n, \olz_1, \ldots, \olz_n  \right)$, according to the bidegrees of homogeneity:
$$
\mcH_k(\mR^{2n};\mC) = \bigoplus_{a+b=k} \ \mcH_{a,b}(\mR^{2n};\mC)
$$
In \cite{FQA} we have established in detail the fundaments of a function theory called {\em quaternionic Clifford analysis} (see also \cite{ABBSS1, ABBSS2, struppa, DES, Ee, PSS}), which is a refinement of Hermitian Clifford analysis (see e.g. \cite{partI, partII, toulouse, eel2, rocha, sabadini}), in its turn a refinement of Euclidean Clifford analysis. Clifford analysis (see e.g. \cite{bds, dss, gilbert, ghs, guerleb}) is, in its most basic form, a generalization to higher dimension of holomorphic function theory in the complex plane. The fundamental group of Euclidean Clifford analysis in $\mR^m$ is the Spin$(m)$ group, which doubly covers the SO$(m)$ group. The fundamental group of Hermitian Clifford analysis in $\mR^{2n}$ is the U$(n)$ group. The corresponding Fischer decompositions in terms of monogenic or Hermitian monogenic homogeneous polynomials respectively, are refinements of the Fischer decompositions  (\ref{harm-decomp}) and (\ref{harmu}) (see also \cite{brno}). As shown in \cite{FQA}, the fundamental group underlying quaternionic Clifford analysis in $\mR^{4p}$ (where the dimension now is a fourfold: $m=2n=4p$), is the symplectic group Sp$( p)$. In order to obtain the corresponding Fischer decomposition it is crucial to know how to further decompose the space  $\mcH_{a,b}(\mR^{2n};\mC)$ as a module for Sp$( p)$. This is the problem we tackle in the present paper.


\section{The symplectic Lie group and Lie algebra}


The symplectic group Sp$( p)$ is the real Lie group of square $p \times p$ matrices with quaternion entries, preserving the symplectic inner product
$$
\langle \xi , \eta \rangle_{\mH} = \xi_1 \overline{\eta}_1 + \xi_2 \overline{\eta}_2 + \cdots +  \xi_p \overline{\eta}_p \quad \xi, \eta \in \mH^p
$$
where $\overline{\cdotp}$ stands for quaternionic conjugation. Equivalently, we can describe Sp$( p)$ as
$$
{\rm Sp}( p) = \{ A \in {\rm GL}_p(\mH) : A A^* = E_p \}
$$
Square matrices in $M_p(\mH)$ may be embedded in $M_{2p}(\mC)$ by the injective homomorphism
$$
\psi_p : M_p(\mH) \longrightarrow M_{2p}(\mC)
$$
where for each quaternion entry
$$
z + w \, j = (x + y \, i) + (u + v \, i) \, j = x + y \, i + u \, j + v \, k
$$
the $2 \times 2$ complex matrix
$\bigl(\begin{smallmatrix}
\phantom{-}z&w\\ - \overline{w}&\olz
\end{smallmatrix} \bigr)$
is substituted.
In this way it turns out that $\psi_p({\rm Sp}( p))$ is a subgroup of SU$(2p)$.\\

At the level of Lie algebra we have the following picture. The real symplectic Lie algebra $\gsp( p)$ of skew--symplectic $M_p(\mH)$ matrices
$$
\gsp( p) = \{ A \in {\rm GL}_p(\mH) : A + A^* = 0  \}
$$
is isomorphic with the subalgebra $\psi(\gsp( p))$ of the Lie algebra $\gu(2p)$ of skew--hermitian $M_{2p}(\mC)$ matrices. Moreover, for $A \in \gsp( p)$, the complex matrix $\psi(A)$ satisfies the relation
\begin{equation}
\label{psip}
\psi(A)^T \, \mI_{2p} + \mI_{2p} \, \psi(A) = 0
\end{equation}
where $\cdotp^T$ stands for the transpose and $\mI_{2p}$ is the complex structure introduced in Section 1.\\[-2mm]

On the other hand, there is the complex symplectic Lie group Sp$_{2p}(\mC)$
of complex linear matrices preserving the standard skew--hermitian form on $\mC^{2p}$:
$$
{\rm Sp}_{2p}(\mC) = \{ A \in {\rm GL}_{2p}(\mC) : A^T \, \mI_{2p} \, A = \mI_{2p}  \}
$$
and its corresponding complex symplectic Lie algebra $\gsp_{2p}(\mC)$ given by
$$
\gsp_{2p}(\mC) = \{ A \in {\rm GL}_{2p}(\mC) : A^T \, \mI_{2p}  +  \mI_{2p} \, A = 0 \}
$$
This Lie algebra  $\gsp_{2p}(\mC)$ is a subalgebra of $\gsl_{2p}(\mC)$; it can be decomposed into the direct sum of its Hermitian subspace and its skew--hermitian subalgebra, both spaces being isomorphic through multiplication by the imaginary unit $i$ :
$$
\gsp_{2p}(\mC) = \left( \gsp_{2p}(\mC)  \cap \gu(2p) \right) \oplus i \left( \gsp_{2p}(\mC)  \cap \gu(2p) \right)
$$
In view of (\ref{psip}) this leads to the following result (see also \cite{FQA}).

\begin{proposition}
The real symplectic Lie algebra $\gsp( p)$ of skew--symplectic $M_p(\mH)$--matrices is isomorphic with the compact form 
 $\gsp_{2p}(\mC)  \cap \gu(2p) $ of the complex symplectic Lie algebra  $\gsp_{2p}(\mC)$:
 $$
 \psi (\gsp( p)) = \gsp_{2p}(\mC)  \cap \gu(2p)
 $$
\end{proposition}
\noindent
Henceforth we will use the Lie algebra $\gsp_{2p}(\mC)$.\\

Now, let us consider the space $\mcH_{a,b}\domcodom$ of complex valued $(a,b)$--homogeneous harmonic polynomials in the variables $(z_1, z_2, \ldots, z_{2p}, \olz_1, \olz_2,\ldots, \olz_{2p})$. Seen the surjectivity of the Laplace operator
$$
\Delta_{4p} : \mcP_{a,b}\domcodom \longrightarrow  \mcP_{a-1,b-1}\domcodom
$$
we have
\begin{eqnarray}
{\rm dim} \, \mcH_{a,b}\domcodom &=& {\rm dim} \, \mcP_{a,b}\domcodom - {\rm dim} \, \mcP_{a-1,b-1}\domcodom\nonumber\\[1mm]
&=& \binom{2p+a-1}{a} \binom{2p+b-1}{b} -  \binom{2p+a-2}{a-1} \binom{2p+b-2}{b-1}\\[1mm]
&=&  \binom{2p+a-1}{2p-1} \binom{2p+b-1}{2p-1} -  \binom{2p+a-2}{2p-1} \binom{2p+b-2}{2p-1}\nonumber
\end{eqnarray}
\noindent
In order to decompose the space $\mcH_{a,b}\domcodom$ into $\gsp_{2p}(\mC)$--irreducibles, use could be made of existing branching rules when restricting $\ggl_{2p}(\mC)$ to $\gsp_{2p}(\mC)$. To that end we have to know the behaviour of $\mcH_{a,b}\domcodom$ as a module for $\ggl_{2p}(\mC)$. Let us recall that in \cite{eel2} the spaces  $\mcH_{a,b}\domcodom$ have been identified as irreducible modules for its simple Lie subalgebra $\gsl_{2p}(\mC)$, with highest weight vector
$$
\mcH_{a,b}\domcodom \cong (a+b, b, \cdots, b)
$$
of length $2p-1$. Interpreted as a representation space for $\ggl_{2p}(\mC)$ we have
$$
\mcH_{a,b}\domcodom \cong (a, 0, \cdots, 0, -b)
$$
instead, where the highest weight vector now has length $2p$. In fact this is telling us that
$$
\mcH_{a,b}\domcodom \cong \big( V^{\boxtimes a}  \big)   \boxtimes \big(\overline{V}^{\boxtimes b}   \big)
$$
where $\boxtimes$ stands for the Cartan product, with $V \cong \mC^{2p}$ the fundamental representation and $\overline{V}$ its dual.
The branching rules when restricting $\ggl_{2p}(\mC)$ to $\gsp_{2p}(\mC)$ could be found in full generality in \cite{HTW}, the branching multiplicities being expressed in terms of Littlewood--Richardson coefficients. However, due to the rather simple highest weight to start with, the actual situation is not that complicated and one obtains, for $a > b$
\begin{equation}
\label{branching}
 \mcH_{a,b} \Biggr\rvert_{\gsp_{2p}(\mC)}^{\ggl_{2p}(\mC)} = (a,b)_s \oplus  (a+1,b-1)_s \oplus \cdots \oplus (a+b-1,1)_s \oplus  (a+b)_s
\end{equation}
where the shorthand notation $(\lambda)_s$ refers to an irreducible representation for $\gsp_{2p}(\mC)$, and stands for a {\em symplectic} highest weight $(\lambda, 0, \ldots, 0)$ of length $p$. Also note that if $a$ or $b$ equals zero then no branching occurs, meaning that in that case 
$ \mcH_{a,b}$ is symplectically irreducible.\\[-2mm]

In order to characterize the spaces of the form $(a,b)_s$ in terms of homogeneous polynomials on $\mR^{4p} \cong \mC^{2p}$, we will establish, in the next section, an alternative realization for $\gsl(2,\mC) \cong {\rm Alg}_{\mC}(X,Y,H)$. To that end we introduce new differential operators appearing in quaternionic Clifford analysis.


\section{Quaternionic Clifford analysis: the basics}


As is well--known, when establishing Hermitian Clifford analysis (see e.g.\ \cite{partI} ) use is made of the projection operators
$$
\onehalf \left(  {\bf 1} \pm i \, \mI_{2n}  \right)
$$
where $\mI_{2n}$ is the standard complex structure on $\mR^{2n}$ (see Section 1), leading to the standard Hermitian vector variables
$$
\uz = \sum_{k=1}^n \ z_k \gf_k \quad {\rm and \quad } \uzd = \sum_{k=1}^n \ \olz_k \gfd_k
$$
and the Hermitian Dirac operators
$$
\upzd = \sum_{k=1}^n \ \p_{\olz_k} \gf_k \quad {\rm and  } \quad \upz = \sum_{k=1}^n \ \p_{z_k} \gfd_k 
$$
where the Witt basis vectors $\left(\gf_k, \gfd_k : k=1,\ldots,n\right)$ are given by
$$
\gf_k = - \onehalf \left(  {\bf 1} - i \, \mI_{2n}  \right)[e_{2k-1}]  \quad {\rm and}  \quad \gfd_k = \onehalf \left(  {\bf 1} + i \, \mI_{2n}  \right)[e_{2k-1}]
$$
$(e_1,\ldots,e_{2n})$ being an orthonormal basis in $\mR^{2n}$.

\begin{lemma} (see \cite{howe})
The Hermitian variables and Dirac operators enjoy the  anti--commutation relations
\begin{eqnarray*}
&& \{ \uz , \uzd \}  =  \nz^2    \\
&& \{ \upz , \upzd \}  =  \frac{1}{4} \, \Delta_{2n} \\
&& \{ \upz , \uz \}   =  \mE_z + \beta  \\
&& \{ \upzd ,  \uzd \}  =  \mE_z^\dagger + n - \beta \\
&& \{ \upz , \uzd \}   = 0  = \{ \upzd ,  \uz \}  
\end{eqnarray*}
where $\beta $ is the so--called spin--Euler operator given by $\beta =  \sum_{k=1}^n \gfd_k \gf_k = n - \sum_{k=1}^n \gf_k \gfd_k$.
They span the odd part of the Lie super algebra $\gsl(1|2) = \gog_0 \oplus \gog_1$ with
\begin{eqnarray*}
\gog_0 & = & \ggl(2,\mC) = \mC \oplus \gsl(2,\mC) = {\rm span}_{\mC}(\mE_z^\dagger - \mE_z + n - 2\beta) \oplus {\rm Alg}_{\mC}(\mE_z^\dagger + \mE_z + n, \onehalf \nz^2, - \onehalf \Delta_{2n}) \\
\gog_1 &=& {\rm span}_{\mC} (\uz, \uzd, \upz, \upzd)
\end{eqnarray*}
\end{lemma}

The central notion in Hermitian Clifford analysis is that of a Hermitian monogenic function, the definition of which is as follows.

\begin{definition}
A differentiable function $F$ defined in a domain $\Omega$ of $\mR^{2n}$ and taking its values in the complex Clifford algebra $\mC_{2n}$ or in spinor space $\mS$, is called {\em Hermitian monogenic} in $\Omega$ if it satisfies the system $\{ \upz F = 0,   \upzd F = 0\}$.
\end{definition}

Now, taking the dimension to be a fourfold: $m=2n=4p$, a {\em quaternionic structure} on $\mR^{4p}$ is established by introducing, next to the standard complex structure $\mI_{4p}$, a second complex structure $\mJ_{4p} \in$ SO$(4p)$ such that $\mJ_{4p}^2 = - E_{4p}$ and $\mI_{4p}$ and $\mJ_{4p}$ are anti--commuting. This second complex structure $\mJ_{4p}$ may be realized as
$$
\mJ_{4p}  = {\rm diag}
\begin{pmatrix}
0 & 0 & 1 & 0\\
0 & 0 & 0 & -1\\
-1 & 0 & 0 & 0 \\
0 & 1 & 0 & 0 
\end{pmatrix}
$$
New twisted Hermitian variables and Dirac operators then are defined by
$$
\uzj = \mJ_{4p} [\uz] = \sum_{k=1}^p \ \gfd_{2k-1} z_{2k} - \gfd_{2k} z_{2k-1}
$$
$$
\uzdj = \mJ_{4p} [\uzd] = \sum_{k=1}^p \ \gf_{2k-1} \olz_{2k} - \gf_{2k} \olz_{2k-1}
$$
$$
\upzj = \mJ_{4p} [\upz] = \sum_{k=1}^p \ \gf_{2k-1} \p_{z_{2k}} - \gf_{2k} \p_{z_{2k-1}}
$$
$$
\upzdj = \mJ_{4p} [\upzd] = \sum_{k=1}^p \ \gfd_{2k-1} \p_{\olz_{2k}} - \gfd_{2k} \p_{\olz_{2k-1}}
$$

\begin{lemma}
The twisted Hermitian variables and Dirac operators enjoy the anti--commutation relations
\begin{eqnarray*}
&& \{ \uzj , \uzdj \}  =  \nz^2    \\
&& \{ \upzj , \upzdj \}  =  \frac{1}{4} \, \Delta_{4p} \\
&& \{ \upzj , \uzj \}   =  \mE_z + 2p - \beta  \\
&& \{ \upzdj ,  \uzdj \}  =  \mE_z^\dagger + \beta \\
&& \{ \upzj , \uzdj \}   =  0  = \{ \upzdj ,  \uzj \}  
\end{eqnarray*}
\end{lemma}

\begin{remark}
Note the similarity of the anti--commutation relations of the twisted Hermitian variables and Dirac operators with those for the standard ones, which, quite naturally, follows from the fact that $\mJ_{4p} \in {\rm SO}(4p)$.
\end{remark}

\begin{remark}
While the operators $\upz$ and $\upzd$ are invariant under {\rm U}$(2p)$, the four Dirac operators $\upz, \upzd, \upzj, \upzdj$, taken together, are invariant under the action of the symplectic group {\rm Sp}$( p)$.
\end{remark}

\begin{definition}
A differentiable function $F : \mR^{4p} \longrightarrow \mS$ is called {\em quaternionic monogenic} in the domain $\Omega \subset \mR^{4p}$ if it is a simultaneous null--solution for the four operators $\upz, \upzd, \upzj, \upzdj$.
\end{definition}

New operators are now arising by considering the mixed anti--commutator relations of the standard and twisted Hermitian variables and Dirac operators. Indeed we can define
\begin{eqnarray*}
&& \SE  :=  \phantom{-} \{  \upzdj  , \uz \} = \sum_{k=1}^p \ z_{2k-1} \p_{\olz_{2k}}  -  z_{2k} \p_{\olz_{2k-1}} \\
&& \SED  :=  - \{ \upzj , \uzd  \} = - \sum_{k=1}^p \ \olz_{2k-1} \p_{z_{2k}}  -  \olz_{2k} \p_{z_{2k-1}}
\end{eqnarray*}
and there also holds
\begin{eqnarray*}
\{\upzd ,  \uzj \} &=& \sum_{k=1}^p \ z_{2k} \p_{\olz_{2k-1}}  -  z_{2k-1} \p_{\olz_{2k}} = - \, \SE \\
\{\upz ,  \uzdj \} &=& \sum_{k=1}^p \ \olz_{2k} \p_{z_{2k-1}}  -  \olz_{2k-1} \p_{z_{2k}} = \phantom{-} \SED
\end{eqnarray*}
These new operators enjoy the following properties.

\begin{lemma}
The operators $\SE$ and $\SED$ are invariant under the symplectic action.
\end{lemma}

\begin{lemma}
One has
$$
\gsl(2,\mC) \cong {\rm Alg}_\mC \left( \mE_z^\dagger - \mE_z , \SED, \SE  \right)
$$
these three generating operators commuting with the harmonic triplet $(H, X, Y)$ introduced in Section 1.
\end{lemma}
\pf
Direct computation shows that indeed:\\[1mm]
(i) $[ \mE_z^\dagger - \mE_z , \SED] = 2 \SED$\\[1mm]
(ii) $[ \mE_z^\dagger - \mE_z, \SE] = - 2 \SE$\\[1mm]
(iii) $[\SED, \SE] =  \mE_z^\dagger - \mE_z$\\[1mm]
(iv) $[\SED, \onehalf \nz^2] = [\SE, \onehalf \nz^2]$ = 0\\[1mm]
(v) $[\SED, - \onehalf \Delta_{4p}] = [\SE, - \onehalf \Delta_{4p}]$ = 0\\[1mm]
As it was already shown that $ \mE_z^\dagger - \mE_z$ commutes with $ \onehalf \nz^2, - \onehalf \Delta_{4p}$ and $\mE_z + \mE_z^\dagger + 2p$, the proof is complete.\eop

\begin{corollary}
\label{so4}
One has
$$
{\rm Alg}_{\mC} \left(   \mE_z + \mE_z^\dagger + 2p, \onehalf   |z|^2, - \onehalf \Delta_{4p} \right) \oplus {\rm Alg}_{\mC}  \left(   \mE_z^\dagger - \mE_z, \SED, \SE \right) \cong \gsl(2,\mC) \oplus \gsl(2,\mC) =  \gso(4,\mC)
$$
\end{corollary}


\section{Symplectic harmonics}


If we consider the operators $\SE$ and $\SED$ as acting between the spaces $\mcH_{a,b}\domcodom$ of  complex valued bi--homogeneous harmonic polynomials, we obtain
$$
\mcH_{0,b+a}  \ \rightleftarrows \ \cdots \ \rightleftarrows \ \mcH_{a-1,b+1} \ \overset{\SE}{\underset{\SED}{\rightleftarrows}} \ \mcH_{a,b} \ \rightleftarrows \ \mcH_{a+1,b-1} \ \rightleftarrows \ \cdots \ \rightleftarrows \ \mcH_{a+b,0} 
$$
and we define the kernel spaces
$$
\HS_{a,b} = \mcH_{a,b}\domcodom  \cap \mbox{Ker} \, \SE  = \mcP_{a,b}\domcodom \cap \mbox{Ker} (\Delta_{4p}, \SE) \quad (a \geq b)
$$
and
$$
\HSD_{a,b} = \mcH_{a,b}\domcodom  \cap \mbox{Ker} \, \SED  = \mcP_{a,b}\domcodom \cap \mbox{Ker} (\Delta_{4p}, \SED) \quad (a \leq b)
$$
These kernel spaces will show to be crucial in the decomposition of $\mcH_{a,b}\domcodom$ in terms of Sp$( p)$--irreducibles. We call their elements {\em (adjoint) symplectic harmonics}. It will be shown further on (see Corollary 4 and Proposition 4) that $\mcH_{a,b}\domcodom \cap \mbox{Ker} \, \SED = \{0\}$ for $a>b$ and $\mcH_{a,b}\domcodom \cap \mbox{Ker} \, \SE = \{0\}$ for $a<b$.

\begin{remark}
For all $k$, $\mcP_{k,0}\domcodom = \mcH_{k,0} = \HS_{k,0}$ and $\mcP_{0,k}\domcodom = \mcH_{0,k} = \HSD_{0,k}$, since the homogeneous polynomials in $\mcH_{k,0}$ (respectively $\mcH_{0,k}$) do not contain the variables $(\olz_1,\ldots,\olz_{2p})$ (respectively $(z_1,\ldots,z_{2p})$).
\end{remark}

With respect to the traditional Fischer inner product, given by
$$
\langle f , g \rangle = f(\upzd,\upz) \ \overline{g} \, \biggr\rvert_{\uz=0}
$$
where $ f(\upzd,\upz)$ is obtained by substituting $\p_{\olz_j}$ for $z_j$ and $\p_{z_j}$ for $\olz_j$ in $f(z_1,\ldots,z_{2p},\olz_1,\ldots,\olz_{2p})$, each of the spaces $\mcH_{a,b}\domcodom$ can be decomposed as the direct sum
\begin{eqnarray*}
\mcH_{a,b} &=& \HS_{a,b} \oplus (\HS_{a,b})^{\perp} \quad a \geq b \\
\mcH_{a,b} &=& \HSD_{a,b} \oplus (\HSD_{a,b})^{\perp} \quad a \leq b
\end{eqnarray*}
where the orthogonal complements  $(\HS_{a,b})^{\perp}$ and  $(\HSD_{a,b})^{\perp}$ are isomorphic with Im$_\SE(\mcH_{a,b})$ and Im$_\SED(\mcH_{a,b})$ respectively. We will now determine those orthogonal complements explicitly.

\begin{lemma}
With respect to the Fischer inner product, the operators $\SE$ and $\SED$ are adjoint operators, i.e. for polynomials $P \in \mcP_{a,b}$ and $Q \in \mcP_{a+1,b-1}$ there holds
$$
\langle \SE P , Q \rangle = \langle P , \SED Q \rangle
$$
\end{lemma}

\pf
It is clear that the Fischer inner product of two monomials is zero unless both monomials are equal up to a constant. This observation and a straightforward calculation lead to the desired result.
\eop

\begin{proposition}
For $a \geq b$, the space $\mcH_{a,b}$ may be decomposed as
\begin{eqnarray}
\label{step1}
\mcH_{a,b} = \HS_{a,b} \oplus \SED \mcH_{a+1,b-1}
\end{eqnarray}
\end{proposition}

\pf
In fact we prove that, with respect to the Fischer inner product,
$
\left( \SED  \mcH_{a+1,b-1} \right)^{\perp} =  \HS_{a,b}.
$
It is important to note that if a function $F$ is harmonic, then also $\SE F$ and $\SED F$ are harmonic, since the Laplace operator commutes with both operators $\SE$ and $\SED$. Let $P \in \HS_{a,b}$, then $\SE P = 0$ and so $\langle P, \SED Q \rangle = 0$ for all $Q \in  \mcH_{a+1,b-1}$, which means that $P$ is orthogonal to $\SED  \mcH_{a+1,b-1}$ or $P \in \left( \SED  \mcH_{a+1,b-1} \right)^\perp$. Conversely, let $P \in \left( \SED  \mcH_{a+1,b-1} \right)^{\perp}$. Then $\SE P \in \mcH_{a+1,b-1}$ and 
$
\langle \SE P , Q \rangle = \langle P , \SED Q \rangle  = 0
$
for all $Q \in \mcH_{a+1,b-1}$. In particular, for $Q = \SE P$ we find
$
\langle \SE P , \SE P \rangle = 0
$
whence $\SE P = 0$ or $P \in \HS_{a,b}$.
\eop

\begin{corollary}
For $a \geq b$, the space $\mcH_{a,b}$ may be decomposed as
\begin{eqnarray}
\label{symplectic-harmonic-decomp}
\mcH_{a,b}  \, = \,  \HS_{a,b}  \, \oplus  \, \SED  \HS_{a+1,b-1}  \, \oplus \SE^{\dagger 2} \HS_{a+2,b-2}  \, \oplus  \, \cdots  \,\oplus  \, \SE^{\dagger b}  \HS_{a+b,0}
\end{eqnarray}
\end{corollary}

\pf
Consecutive application of the decomposition (\ref{step1}) leads to the desired result.
\eop

\begin{lemma}
One has for $H_{\alpha, \beta} \in \mcH_{\alpha, \beta}$\\
$$\SE \SE^{\dagger k} H_{\alpha,\beta} = k(\alpha - \beta - k +1) \SE^{\dagger (k-1)} H_{\alpha, \beta} + \SE^{\dagger k} \SE H_{\alpha,\beta}$$
and in particular for $H_{\alpha, \beta}^S \in \HS_{\alpha, \beta}$
$$\SE \SE^{\dagger k} H_{\alpha,\beta}^S = k(\alpha - \beta - k +1) \SE^{\dagger (k-1)} H_{\alpha, \beta}^S$$
and $$\SE^\ell \SE^{\dagger k} H_{\alpha,\beta}^S = k(k-1) \cdots (k-\ell+1)(\alpha - \beta - k +1)\cdots(\alpha - \beta - k +\ell) \SE^{\dagger (k-\ell)} H_{\alpha, \beta}^S$$
\end{lemma}

\pf
Straightforward computation based on the commutator $[\SE , \SED] = \mE_z - \mE_z^\dagger$ (see Lemma 4).
\eop

\begin{corollary}
For $a \geq b$, the mappings
$$ 
\begin{array}{clllr}
\SE &:& \SED \HS_{a+1,b-1} &\longrightarrow& \HS_{a+1,b-1}\\[2mm]
\SE &:& \SE^{\dagger 2} \HS_{a+2,b-2} &\longrightarrow& \SED \HS_{a+2,b-2}\\
\vdots &&\\[1mm]
\SE &:& \SE^{\dagger b} \HS_{a+b,0} &\longrightarrow& \SE^{\dagger (b-1)} \HS_{a+b,0}
\end{array}
$$
are isomorphisms, their inverses being, up to constants, restrictions of the operator $\SED$ to the corresponding spaces.
\end{corollary}

\pf
We prove that for $j=1,\ldots,b$
$$
\SE \, : \, \SE^{\dagger j} \HS_{a+j,b-j} \longrightarrow \SE^{\dagger (j-1)} \HS_{a+j,b-j}\\
$$
is an isomorphism. First take $g \in \SE^{\dagger (j-1)} \HS_{a+j,b-j}$, meaning that $g =  \SE^{\dagger (j-1)} h$ with $h \in \HS_{a+j,b-j}$, and consider 
$$
f = \frac{1}{j(a-b+j+1)} \, \SED g = \frac{1}{j(a-b+j+1)} \, \SE^{\dagger j} h  \in \SE^{\dagger j} \HS_{a+j,b-j}
$$ 
Then, using the formulae of Lemma 6, it follows that $\SE f = \SE^{\dagger (j-1)} h = g$, and so the considered mapping is surjective. Moreover $\left(\SE^{\dagger j} \HS_{a+j,b-j}\right)^\perp = \HS_{a,b} \oplus \SED \HS_{a+1,b-1} \oplus \cdots \oplus \SE^{\dagger (j-1)} \HS_{a+j-1,b-j+1}$ implying that this mapping is also injective. Clearly
$$
\left(  \SE\biggr\rvert_{(\SE^{\dagger j} \HS_{a+j,b-j})}  \right)^{-1}  = \frac{1}{j(a-b+j+1)} \ \SED
$$
\eop

\begin{corollary}
For $a \geq b$ and $j = 1, \ldots, b$ one has
$$\mcH_{a+j,b-j}  \cap  {\rm Ker} \, \SED = \{0\}$$
\end{corollary}

\pf
Take $f \in \mcH_{a+j,b-j}$ with $\SED f = 0$ and hence also $\SE \SED f = 0$. In view of Corollary 2, the function $f$ can be decomposed as $f = \sum_{k=0}^{b-j} \, f_k$ with $f_k \in \SE^{\dagger k} \HS_{a+j+k,b-j-k}$. It then follows, in view of Corollary 3, that $\SE \SED f = 0$ implies $f_0 = f_1 = \ldots = f_{b-j} = 0$, and hence also $f = 0$.
\eop\\

In a similar way as for the case where $a \geq b$, the following results hold for the case where $a \leq b$.

\begin{proposition}
\label{S-decomp-dagger}
For $a \leq b$ one has
\begin{itemize}

\item[(i)] 
the space $\mcH_{a,b}$ may be decomposed as
\begin{eqnarray}
\label{symplectic-harmonic-decomp-bis}
\mcH_{a,b} & = & \HSD_{a,b} \oplus \SE \mcH_{a-1,b+1} \nonumber \\
 & = & \HSD_{a,b} \oplus \SE  \HSD_{a-1,b+1} \oplus \SE^{ 2} \HSD_{a-2,b+2} \oplus \cdots \oplus \SE^{a}  \HSD_{0,b+a}
\end{eqnarray}

\item[(ii)]
the mappings
$$
\begin{array}{ccccr}
\SED &:& \SE \HSD_{a-1,b+1} &\longrightarrow& \HSD_{a-1,b+1}\\
\SED &:& \SE^{ 2} \HSD_{a-2,b+2} &\longrightarrow& \SE \HSD_{a-2,b+2}\\
\vdots &&\\
\SED &:& \SE^{a} \HSD_{0,b+a} &\longrightarrow& \SE^{(a-1)} \HSD_{0,b+a}
\end{array}
$$
are isomorphisms, their inverses being, up to constants, restrictions of the operator $\SE$ to the corresponding spaces.

\item[(iii)]
$\mcH_{a-j,b+j}  \cap  {\rm Ker} \, \SE = \{0\}$ for $j=1,\ldots,a$.

\end{itemize}
\end{proposition}

\begin{corollary}
For $a\geq b$ there holds
\begin{eqnarray*}
{\rm dim} \ \HS_{a,b}\domcodom & = & {\rm dim} \ \mcH_{a,b}  - {\rm dim}  \ \mcH_{a+1,b-1}\\
& = &  {\rm dim}  \ \mcP_{a,b}  - {\rm dim} \  \mcP_{a-1,b-1} -  {\rm dim} \ \mcP_{a+1,b-1} +  {\rm dim} \ \mcP_{a,b-2}\\
& = & \frac{(2p-1)(2p-2)(a-b+1)(a+b+2p-1)(a+2p-2)!(b+2p-3)!}{((2p-1)!)^2(a+1)!b!}
\end{eqnarray*}
For $a\leq b$ there holds
\begin{eqnarray*}
{\rm dim} \ \HSD_{a,b\domcodom } & = & {\rm dim} \ \mcH_{a,b}  - {\rm dim}  \ \mcH_{a-1,b+1}\\
& = &  {\rm dim}  \ \mcP_{a,b}  - {\rm dim} \  \mcP_{a-1,b-1} -  {\rm dim} \ \mcP_{a-1,b+1} +  {\rm dim} \ \mcP_{a-2,b}\\
& = & {\rm dim} \ \HS_{b,a}\domcodom
\end{eqnarray*}
\end{corollary}

\noindent As ${\rm dim} \, \HSD_{a,b} = {\rm dim} \, \HS_{b,a}$, the spaces  $\HSD_{a,b}$ and  $\HS_{b,a}$ are isomorphic. Obviously this isomorphism is realized by complex conjugation which, indeed, maps  $\mcH_{a,b}$ and  $\mcH_{b,a}$ onto each other, since the Laplace operator is invariant under complex conjugation, and moreover the operators $\SE$ and $\SED$ are complex conjugated up to a minus sign. \\[-2mm]

There is, however, another -nice- way to express this isomorphism, which is closely related to the quaternionic structure  $(\mI, \mJ, \mK)$ introduced in \cite{FQA} to study the fundaments of quaternionic Clifford analysis (see Section 3). For a function $F(z_1,\ldots,z_{2p},\olz_1,\ldots,\olz_{2p})$ consider the transformation $T$, mapping $F$ onto the function $T[F]$ by substituting for the variables $z_{2k-1}, z_{2k},\olz_{2k-1},\olz_{2k}$ the variables $-\olz_{2k}, \olz_{2k-1},-z_{2k},z_{2k-1} (k=1,\ldots,p)$ respectively. In fact this is the transformation associated to the second complex structure $\mJ_{4p} \in {\rm SO}(4p)$ in the quaternionic structure. If $h_{a,b} \in \mcH_{a,b}$, then $T[h_{a,b}] \in  \mcH_{b,a}$ since $T$ commutes with the Laplace operator. Let us now compute the commutation relations of $T$ with the operators $\SE$ and $\SED$.

\begin{lemma}
For the  transformation $T$ introduced above, it holds
$$
\SED \, T = - T \, \SE \quad {\rm and} \quad \SE \, T = - T \SED
$$
\end{lemma}

\pf
We consecutively have 
\begin{eqnarray*}
\SED T[F] & = & \sum_{k=1}^p (\olz_{2k} \p_{z_{2k-1}}  -  \olz_{2k-1} \p_{z_{2k}} ) T[F]\\
& = & \sum_{k=1}^p \olz_{2k} T[\p_{\olz_{2k}}F] -  \olz_{2k-1} T[-\p_{\olz_{2k-1}}F]\\
& = & T\left[\sum_{k=1}^p - z_{2k-1} \p_{\olz_{2k}}F +  z_{2k} \p_{\olz_{2k-1}}F\right]\\
& = & T[ - \SE F]
\end{eqnarray*}
Next, taking into account that $T^2 = - {\bf 1}$, we also have
$$
T \, \SED \, T \, T = - T \, T \, \SE \, T \quad {\rm or} \quad T \, \SED = - \SE \, T
$$
\eop

\begin{corollary}
\label{iso}
If $F$ is in {\em Ker}\,$\SE$, then $T[F]$ is in {\em Ker}\,$\SED$ and vice versa, and, consequently
$$ 
T : \HS_{a,b} \longleftrightarrow \HSD_{b,a}
$$
is an isomorphism. 
\end{corollary}

\begin{remark}
Taking for the operator $T$ the operator associated with the third complex structure $\mK_{4p}$, which corresponds to the change of variables
$
z_{2k-1} \mapsto i \olz_{2k}  ,  z_{2k} \mapsto - i \olz_{2k-1} , \olz_{2k-1} \mapsto  - i z_{2k} , \olz_{2k} \mapsto i z_{2k-1}, 
$
 also leads to an isomorphism between the spaces $\HS_{a,b}$ and $\HSD_{b,a}$. The operator associated to the first complex structure $\mI_{4p}$ is an automorphism of both spaces $\HS_{a,b}$ and $\HSD_{a,b}$.
\end {remark}

In the special case where $a=b=k$, the isomorphism between $\HS_{k,k}$ and $\HSD_{k,k}$ becomes the identity. This is a special case (for $j=0$) of the following lemma, the proof of which invokes the Fischer decomposition established in the next section (see Theorems \ref{theofisch1}--\ref{theofisch2}).

\begin{lemma}
\label{identity}
For all $j=0,1,\ldots,k$ one has
$$
(\SED)^j \HS_{k+j,k-j} = (\SE)^j \HSD_{k-j,k+j}
$$
\end{lemma}

\pf
Since $(\SED)^j \HS_{k+j,k-j}$ and  $(\SE)^j \HSD_{k-j,k+j}$ both are non--trivial, irreducible $\gsp$-submodules of $\mcH_{k,k}$ with the same highest weight $(k+j,k-j)_s$, they coincide seen the Fischer decomposition of $\mcH_{k,k}$.
\eop \\[-2mm]

\noindent Also the case where $a-b=1$ is interesting, and is obtained (by taking $j=0$) from the following lemma, which also leans upon the Fische decomposition.

\begin{lemma}
\label{isomorphism}
For all $j=0,1,\ldots,k$ one has
$$
(\SED)^{j+1} \HS_{k+1+j,k-j} = (\SE)^j \HSD_{k-j,k+1+j}
$$
\end{lemma}

\pf
First note that $(\SED)^{j+1} \HS_{k+1+j,k-j}$ is not the null--space, since $(\SED)^{j} \HS_{k+1+j,k-j} \neq 0$ and\\
Ker $\SED \cap \mcH_{k+1,k} = 0$. Similarly, also $(\SE)^j \HSD_{k-j,k+1+j}$ is not the null--space. Since both spaces 
$(\SED)^{j+1} \HS_{k+1+j,k-j}$ and  $(\SE)^j \HSD_{k-j,k+1+j}$ are, non--trivial, irreducible $\gsp$--submodules of $\mcH_{k,k+1}$  with the same highest weight $(k+1+j,k-j)$, they coincide seen the Fischer decomposition of $\mcH_{k,k+1}$. \eop\\[-2mm]

\begin{corollary}
For $a > b$, the mappings
$$
\SE^{a-b} : \HSD_{b,a} \longrightarrow \HS_{a,b}
$$
and
$$
\SE^{\dagger (a-b)} : \HS_{a,b} \longrightarrow \HSD_{b,a}
$$
are isomorphisms.
\end{corollary}


\section{Fischer decompositions}


First assume that $a>b$ and compare the decomposition (\ref{symplectic-harmonic-decomp}) for $\mcH_{a,b}$ in terms of symplectic harmonics, viz.
$$
\mcH_{a,b} = \HS_{a,b} \oplus \SED  \HS_{a+1,b-1} \oplus \SE^{\dagger 2} \HS_{a+2,b-2} \oplus \cdots \oplus \SE^{\dagger b}  \HS_{a+b,0}
$$
with the branching (\ref{branching}):
$$
 \mcH_{a,b} \Biggr\rvert_{\gsp_{2p}(\mC)}^{\ggl_{2p}(\mC)} = (a,b)_s \oplus  (a+1,b-1)_s \oplus \cdots \oplus (a+b-1,1)_s \oplus  (a+b)_s
$$
It is then rather straightforward to conjecture that for $a>b$
$$
 (a,b)_s \cong \HS_{a,b}(\mR^{4p};\mC)
$$
That this indeed is the case is shown in the next theorem.

\begin{theorem}
One has, with $a \geq b$,
$$
 (a,b)_s \cong \HS_{a,b}(\mR^{4p};\mC)
$$
where  $(a,b)_s = (a,b,0,\ldots,0)_s$ stands for an irreducible $\gsp_{2p}(\mC)$--representation, the highest weight being of length $p$, and 
$$
\mcH_{a,b}(\mR^{4p};\mC) = \HS_{a,b} \oplus \SED  \HS_{a+1,b-1} \oplus \SE^{\dagger 2}  \HS_{a+2,b-2} \oplus \cdots \oplus \SE^{\dagger b}  \HS_{a+b,0} 
$$
is the Fischer decomposition of the space of complex valued bi--homogeneous harmonic  polynomials in terms of $\gsp_{2p}(\mC)$--irreducibles of complex valued bi--homogeneous symplectic harmonic  polynomials.
\label{theofisch1}
\end{theorem}

\pf
We proceed by induction on $b$.\\
For $b=0$ the result is trivial; indeed, as was already noticed in Section 2, in this case no branching occurs and $\mcH_{a,0}$ is symplectically irreducible.\\
Assuming that the theorem is true for $b-1$ means that
$$
 (a+1,b-1)_s \cong \HS_{a+1,b-1}(\mR^{4p};\mC)
$$
and that 
$$
\mcH_{a+1,b-1}(\mR^{4p};\mC) = \HS_{a+1,b-1} \oplus \SED  \HS_{a+2,b-2} \oplus \SE^{\dagger 2}  \HS_{a+3,b-3} \oplus \cdots \oplus \SE^{\dagger (b-1)}  \HS_{a+b,0} 
$$
is an $\gsp_{2p}(\mC)$--irreducible decomposition, as then is also the case  for
$$
\SED \mcH_{a+1,b-1}(\mR^{4p};\mC) = \SED \HS_{a+1,b-1} \oplus \SE^{\dagger 2}  \HS_{a+2,b-2} \oplus \SE^{\dagger 3}  \HS_{a+3,b-3} \oplus \cdots \oplus \SE^{\dagger b}  \HS_{a+b,0} 
$$
which in fact also reads
$$
\SED \mcH_{a+1,b-1} \Biggr\rvert_{\gsp_{2p}(\mC)}^{\ggl_{2p}(\mC)} = (a+1,b-1)_s \oplus  (a+2,b-2)_s \oplus \cdots \oplus (a+b-1,1)_s \oplus  (a+b)_s
$$
In view of  the branching (\ref{branching}) and the decomposition (\ref{step1}) it follows that
$$
 (a,b)_s \cong \HS_{a,b}(\mR^{4p};\mC)
$$
which finishes the proof. \eop\\[-2mm]

If $a<b$ then we have to compare the decomposition  (\ref{symplectic-harmonic-decomp-bis}) 
$$
\mcH_{a,b} = \HSD_{a,b} \oplus \SE  \HSD_{a-1,b+1} \oplus \SE^2  \HSD_{a-2,b+2} \oplus \cdots \oplus \SE^{a}  \HSD_{0,b+a}
$$
with the branching rule
$$
 \mcH_{a,b} \Biggr\rvert_{\gsp_{2p}(\mC)}^{\ggl_{2p}(\mC)} = (b+a)_s \oplus  (b+a-1,1)_s \oplus \cdots \oplus   (b,a)_s
$$
leading to the complementary conjecture for $a<b$:
$$
 (b,a)_s \cong \HSD_{a,b}(\mR^{4p};\mC)
$$
which is proven in a similar way.

\begin{theorem}
One has, with $a \leq b$,
$$
 (b,a)_s \cong \HSD_{a,b}(\mR^{4p};\mC)
$$
and
$$
\mcH_{a,b}(\mR^{4p};\mC) = \HSD_{a,b} \oplus \SE  \HSD_{a-1,b+1} \oplus \SE^2  \HSD_{a-2,b+2} \oplus \cdots \oplus \SE^{a}  \HSD_{0,b+a} 
$$
is the Fischer decompositions of the space of complex valued bi--homogeneous harmonic  polynomials in terms of $\gsp_{2p}(\mC)$--irreducibles of complex valued bi--homogeneous adjoint symplectic harmonic  polynomials.
\label{theofisch2}
\end{theorem}

\begin{corollary}
With $a>b$, the spaces $\HS_{a,b}(\mR^{4p};\mC)$ and $\HSD_{b,a}(\mR^{4p};\mC)$ are isomorphic irreducible representations for $\gsp_{2p}(\mC)$.
\end{corollary}

We already know the dimension of the spaces $\HS_{a,b}$ and $\HSD_{b,a}$ to be (see Corollary 5)
$$
\mbox{dim} \ \HS_{a,b} = \mbox{dim} \ \HSD_{b,a} = \frac{(2p-1)(2p-2)(a-b+1)(a+b+2p-1)(a+2p-2)!(b+2p-3)!}{((2p-1)!)^2(a+1)!b!}
$$
Now we are able to calculate this dimension in the following alternative way. If $\Gamma_\lambda$ denotes an irreducible representation for $\gsp_{2p}(\mC)$ with highest weight $\lambda = (\lambda_1 \geq \lambda_2  \geq  \cdots  \geq  \lambda_p)$, then (see \cite{fulton}, p. 406)
$$
\mbox{dim} \ \Gamma_\lambda = \prod_{i<j} \ \frac{\ell_i^2-\ell_j^2}{m_i^2-m_j^2} \prod_i \frac{\ell_i}{m_i}
$$
with $\ell_i = \lambda_i + m_i$ and $m_i = p-i+1, i=1,\ldots,p$.
For the highest weight $(a,b)_s$ we have
$$
m(m1,\ldots,m_p) = (p, p-1, \ldots, 1)
$$
and
$$
\ell(\ell_1,\ldots,\ell_p) = (a+p, b+p-1, p-2, p-3, \ldots, 1)
$$
A straightforward calculation then leads to
$$
\mbox{dim} \ (a,b)_s = \frac{(a-b+1)(a+b+2p-1)(a+2p-2)!(b+2p-3)!}{(2p-3)!(a+1)!b!}
$$
which is, quite naturally, the dimension of $\HS_{a,b}$ and $\HSD_{b,a}$.\\

Let us give an illustrative example of the Fischer decompositions above. Take $p=2$, and consider the decompositions
\begin{eqnarray*}
\mcH_{2,2}(\mR^{8};\mC) &=& \HS_{2,2} \oplus \SED  \HS_{3,1} \oplus \SE^{\dagger 2}  \HS_{4,0} \\
\mcH_{2,2}(\mR^{8};\mC) &=& \HSD_{2,2} \oplus \SE  \HSD_{1,3} \oplus \SE^2  \HSD_{0,4}
\end{eqnarray*}
The harmonic polynomial $z_3^2 \olz_1^2 \in \mcH_{2,2}$ is decomposed as
$$
z_3^2 \olz_1^2 = P_1 + P_2 +P_3
$$
with
\begin{eqnarray*}
P_1 &=& \frac{1}{3} z_3^2 \olz_1^2 + \frac{1}{3} z_2^2 \olz_4^2 + \frac{2}{3} z_2 z_3 \olz_1 \olz_4\\
P_2 &=& \onehalf z_3^2 \olz_1^2 - \onehalf z_2^2 \olz_4^2\\
P_3 &=& \frac{1}{6} z_3^2 \olz_1^2 + \frac{1}{6} z_2^2 \olz_4^2 - \frac{2}{3} z_2 z_3 \olz_1 \olz_4
\end{eqnarray*}
The polynomial $P_1$ belongs to $\HS_{2,2} \equiv \HSD_{2,2}$. The polynomial $P_2$ can be written as either
$$
P_2 = \SED Q_2 \quad {\rm with} \quad Q_2 =  \onehalf  (- z_2 z_3^2 \olz_1 - z_2^2 z_3 \olz_4)
$$
or
$$
P_2 = \SE \widetilde{Q_2} \quad {\rm with} \quad  \widetilde{Q_2} = \onehalf ( z_3 \olz_1^2 \olz_4 + z_2 \olz_1 \olz_4^2  )
$$
The polynomial $Q_2$ belongs to $\HS_{3,1}$, while the polynomial $\widetilde{Q_2}$ belongs to $\HSD_{1,3}$.
The polynomial $P_3$ can be written as either
$$
P_3 = \SE^{\dagger 2} Q_3 \quad {\rm with} \quad Q_3 = \frac{1}{12} z_2^2 z_3^2
$$
or
$$
P_3 = \SE^2 \widetilde{Q_3} \quad {\rm with} \quad  \widetilde{Q_3} = \frac{1}{12} \olz_1^2 \olz_4^2
$$
The polynomial $Q_3$ belongs to $\HS_{4,0} \equiv \mcH_{4,0}$, while the polynomial $\widetilde{Q_3}$ belongs to $\HSD_{0,4} \equiv \mcH_{0,4}$.

\begin{corollary}
The Fischer decomposition of the space $\mcP_{a,b}(\mR^{4p};\mC)$ of complex valued bi--homogeneous  polynomials in terms of irreducible symplectic modules, is given by
\begin{equation}
\label{symplectic-decomp-one}
\mcP_{a,b} = \bigoplus_{j=0}^b \ |\uz|^{2j} \ \mcH_{a-j,b-j} =
\bigoplus_{j=0}^b \ \bigoplus_{t=0}^{b-j} \  |\uz|^{2j} \ \SE^{\dagger t} \HS_{a-j+t,b-j-t} \quad (a \geq b)
\end{equation}
or
\begin{equation}
\label{symplectic-decomp-two}
\mcP_{a,b} = \bigoplus_{j=0}^a \ |\uz|^{2j} \ \mcH_{a-j,b-j} =
\bigoplus_{j=0}^a \ \bigoplus_{t=0}^{a-j} \  |\uz|^{2j} \ \SE^{t} \HSD_{a-j-t,b-j+t} \quad (a \leq b)
\end{equation}
\end{corollary}

\begin{corollary}
The space $\mcP(\mR^{4p};\mC)$ may be decomposed in terms of irreducible symplectic modules according to the following diagrams.\\

\noindent
For $\mcP_0(\mR^{4p};\mC)$:
$$
\begin{array}{ccc}
& \HS_{0,0} &
\end{array}
$$

\noindent
For $\mcP_2(\mR^{4p};\mC)$:
$$
\begin{array}{ccc}
& r^2 \HS_{0,0} & \\ \\
\HSD_{0,2} & \SED \HS_{2,0} & \HS_{2,0}\\
& \HS_{1,1} &
\end{array}
$$

\noindent
For $\mcP_4(\mR^{4p};\mC)$:
$$
\begin{array}{ccccc}
&& r^4 \HS_{0,0} && \\ \\
& r^2 \HSD_{0,2} & r^2 \SED \HS_{2,0} & r^2 \HS_{2,0} &\\
&& r^2\HS_{1,1} &&\\ \\
\HSD_{0,4} & \SE \HSD_{0,4} & \SE^{\dagger 2}  \HS_{4,0} & \SED  \HS_{4,0} &  \HS_{4,0}\\
& \HSD_{1,3} & \SED \HS_{3,1} &  \HS_{3,1} &\\
&&  \HS_{2,2} &&
\end{array}
$$

\noindent
etc. for even degree polynomials.\\

\noindent
For $\mcP_1(\mR^{4p};\mC)$:
$$
\begin{array}{cccc}
& \HSD_{0,1} & \HS_{1,0} &
\end{array}
$$

\noindent
For $\mcP_3(\mR^{4p};\mC)$:
$$
\begin{array}{cccc}
& r^2 \HSD_{0,1} & r^2 \HS_{1,0} &\\ \\
\HSD_{0,3} & \SE \HSD_{0,3} & \SED \HS_{3,0} & \HS_{3,0}\\
& \HSD_{1,2} & \HS_{2,1} &
\end{array}
$$

\noindent
For $\mcP_5(\mR^{4p};\mC)$:
$$
\begin{array}{cccccc}
&& r^4 \HSD_{0,1} & r^4 \HS_{1,0} &&\\ \\
& r^2 \HSD_{0,3} & r^2 \SE \HSD_{0,3} & r^2 \SED \HS_{3,0} & r^2 \HS_{3,0} &\\
&& r^2 \HSD_{1,2} & r^2 \HS_{2,1} &&\\ \\
\HSD_{0,5} & \SE \HSD_{0,5}  & \SE^{ 2} \HSD_{0,5} & \SE^{\dagger 2} \HS_{5,0} & \SED \HS_{5,0} & \HS_{5,0}\\
& \HSD_{1,4} & \SE \HSD_{1,4}  &  \SED \HS_{4,1} & \HS_{4,1} &\\
&& \HSD_{2,3} & \HS_{3,2} &
\end{array}
$$

\noindent
etc. for odd degree polynomials.
\end{corollary}


\section{Howe dual pair}


The Fischer decomposition (\ref{harm-decomp}) of the space $\mcP(\mR^m;\mC)$ of complex valued polynomials in terms of spherical harmonics, viz.
$$
\mcP(\mR^m;\mC) = \bigoplus_{k=0}^{\infty} \bigoplus_{p=0}^{\infty} r^{2p} \ \mcH_k(\mR^m;\mC)
$$
shows the drawback that it is not multiplicity-free: each of the SO($m$)--irreducible invariant subspaces $\mcH_k(\mR^m;\mC)$ appears with an infinite multiplicity, since all of 
$$r^{2p} \ \mcH_k \quad , \quad p \in \mN_0$$
$k \in \mN_0$ being fixed, are isomorphic as SO($m$)--modules. Expressing irreducibility with respect to 
$\gog \times {\rm SO}(m)$, $\gog$ being an appropriate Lie algebra, aims at collecting the infinitely many copies of $\mcH_k$ into one single irreducible representation. The so--called Howe dual pair $\left( \mbox{SO}(m),\gog \right)$ is to be found with respect to a bigger algebra in which $\gso(m)$ and $\gog$ are commuting. Seen the action (\ref{action}) of the operators  $X := \frac{1}{2} \ r^2$, $Y := - \frac{1}{2} \ \Delta_m$ and $H := \mE + \frac{m}{2}$, the Lie algebra $\gog$ in this case is $\gsl(2,\mC)$. More background information is to be found in \cite{howe}.\\

Similarly, the Fischer decomposition (\ref{harmu}) of the space $\mcP(\mR^{2n};\mC)$ in terms of Hermitian spherical harmonics:
\begin{equation}
\mcP(\mR^{2n};\mC) =
\bigoplus_{k=0}^{\infty} \bigoplus_{p=0}^{\infty}   \bigoplus_{a=0}^{k} \ r^{2p} \ \mcH_{a,k-a}(\mR^{2n};\mC)
\end{equation}
is not multiplicity free since, for all $a$ and $b$,
$$
r^{2p} \ \mcH_{a,b}, \quad p=0,1,2,\ldots
$$
are isomorphic as U$(n)$--modules. It turns out that the Howe dual pair here is $(\mbox{U}(n), \ggl(2,\mC))$ (see also \cite{howe}), with 
\begin{eqnarray*}
\ggl(2,\mC) & = & \gsl(2,\mC) \oplus \mC\\
\phantom{\ggl(2,\mC)} & = & {\rm Alg}_{\mC}(H,X,Y) \oplus {\rm span}_{\mC}(\mE_z^{\dagger} - \mE_z + n)
\end{eqnarray*}

Now let us have a look at the symplectic Fischer decomposition of the space $\mcP(\mR^{4p};\mC)$. By means of (\ref{symplectic-decomp-one}) and (\ref{symplectic-decomp-two})  we obtain

\begin{equation*}
\mcP(\mR^{4p};\mC)  = \bigoplus_{t=0}^{\infty} \quad
\bigoplus_{a<b} \ |\uz|^{2t} 
\left(  \HSD_{a,b}    \oplus \SE \HSD_{a-1,b+1} \oplus \cdots \oplus \SE^a \HSD_{0,b+a}     
\right)
\end{equation*}
\begin{equation} 
\label{symplectic-decomp}    
\oplus \ \
\bigoplus_{a \geq b} \ |\uz|^{2t} 
\left(  \HS_{a,b}    \oplus \SED \HS_{a+1,b-1} \oplus \cdots \oplus \SE^{\dagger b} \HS_{a+b,0}    
\right)                                  
\end{equation}

or, alternatively
\begin{equation}
\label{symplectic-decomp-bis}
\mcP(\mR^{4p};\mC)  = \bigoplus_{t=0}^{\infty} \quad
\bigoplus_{a \geq b} \quad
\bigoplus_{s=0}^{a-b} \ |\uz|^{2t}  \, \SE^{\dagger s}   \HS_{a,b}  
\end{equation}

or still, interchanging the role of the operators $\SE$ and $\SED$,
\begin{equation}
\label{symplectic-decomp-ter}
\mcP(\mR^{4p};\mC)  = \bigoplus_{t=0}^{\infty} \quad
\bigoplus_{a \leq b} \quad
\bigoplus_{s=0}^{b-a} \ |\uz|^{2t}  \, \SE^s   \HSD_{a,b}  
\end{equation}

Clearly these decompositions are not multiplicity free. Assuming $a > b$, the isomorphic Sp$( p)$--modules may be gathered in the following way
$$
\hspace*{-5mm}
\begin{array}{cccccccccc}
\vdots && \vdots && \vdots & \vdots && \vdots && \vdots\\
\uparrow && \uparrow && \uparrow & \uparrow && \uparrow && \uparrow\\
|\uz|^4 \HS_{a,b} & \rightarrow & |\uz|^4 \SED \HS_{a,b} & \rightarrow   \cdots  \rightarrow & |\uz|^4 \SE^{\dagger \alpha} \HS_{a,b}
\overset{iso}{ \rightarrow } &
|\uz|^4 \SE^\alpha \HSD_{b,a} & \rightarrow \cdots  \rightarrow & |\uz|^4 \SE \HSD_{b,a} & \rightarrow & |\uz|^4 \HSD_{b,a}\\
\uparrow && \uparrow && \uparrow & \uparrow && \uparrow && \uparrow\\
|\uz|^2 \HS_{a,b} & \rightarrow & |\uz|^2 \SED \HS_{a,b} & \rightarrow  \cdots  \rightarrow & |\uz|^2 \SE^{\dagger \alpha} \HS_{a,b}
\stackrel{iso}{ \rightarrow } &
|\uz|^2 \SE^\alpha \HSD_{b,a} & \rightarrow \cdots  \rightarrow & |\uz|^2 \SE \HSD_{b,a} & \rightarrow &|\uz|^2 \HSD_{b,a}\\
\uparrow && \uparrow && \uparrow & \uparrow && \uparrow && \uparrow\\
\phantom{|z|^2} \HS_{a,b} & \rightarrow & \phantom{|z|^2}\SED \HS_{a,b}  & \rightarrow  \cdots  \rightarrow & \phantom{|z|^2} \SE^{\dagger \alpha} \HS_{a,b}
\stackrel{iso}{ \rightarrow } &
\phantom{|\uz|^2} \SE^\alpha \HSD_{b,a} & \rightarrow \cdots  \rightarrow & \phantom{|\uz|^2} \SE \HSD_{b,a} & \rightarrow & \phantom{|\uz|^2} \HSD_{b,a}
\end{array}
$$
where $\alpha = \lfloor \frac{a-b}{2} \rfloor$, and $iso$ is either the identity if $a-b$ is even (see Lemma \ref{identity}), or the mapping $\SED$ if $a-b$ is odd (see Lemma \ref{isomorphism}). In this scheme the horizontal arrows represent the action of the operator $\SED$, while the vertical arrows correspond to multiplication by $|\uz|^2$. If $a<b$ this scheme has to be reinterpreted ''from right to left'', the horizontal arrows, now oriented from right to left, then corresponding to the action of the operator $\SE$. In the special case where $a=b$, the scheme reduces to

$$
\begin{array}{ccc}
&\vdots&\\
|\uz|^4 \HS_{a,a} & = & |\uz|^4 \HSD_{a,a}\\
&\uparrow&\\
|\uz|^2 \HS_{a,a} & = &|\uz|^2 \HSD_{a,a}\\
&\uparrow&\\
\phantom{|\uz|^2} \HS_{a,a} & = & \phantom{|\uz|^2} \HSD_{a,a}\\
\end{array}
$$

Apparently the Howe dual $\gog$ is generated by the operators
$$
\SE , \SED , \mE_z^\dagger - \mE_z , |\uz|^2 , \Delta_{4p} ,  \mE_z + \mE_z^\dagger + 2p
$$
which leads to (see also Corollary \ref{so4})
\begin{eqnarray*}
\gog & = & {\rm Alg}_{\mC} \left(   \mE_z + \mE_z^\dagger + 2p, \onehalf   |\uz|^2, - \onehalf \Delta_{4p} \right) \oplus {\rm Alg}_{\mC}  \left(   \mE_z^\dagger - \mE_z, \SED, \SE \right) \\
& = & \gsl(2,\mC) \oplus \gsl(2,\mC)\\
& = & \gso(4,\mC)
\end{eqnarray*}

So let us decompose the Sp$( p)$--module $\mcP(\mR^{4p};\mC)$ under the combined action of the Howe dual pair $(\gsl(2,\mC) \oplus  \gsl(2,\mC) ) \times {\rm Sp}( p)$. For each irreducible Sp$( p)$--module $\HS_{a,b}$ we choose a basis $\{ H^S_{a,b;j} : j=1,2,\ldots, \mbox{dim}\HS_{a,b} \}$; this is a set of singular vectors, labeled by three parameters $a, b$ and $j$. The repeated action of $X = \onehalf |\uz|^2$ then generates the module $\mV_{a,b;j}$ given by
$$
\mV_{a,b;j} = {\rm span}_{\mC} \{  X^t H^S_{a,b;j} : t = 0,1,2, \ldots    \}
$$
Each of the spaces $\mV_{a,b;j}$ is a realization of a so--called Verma module, an infinite dimensional irreducible $\gsl(2,\mC)$--module, which we denote by $\mI^\infty_{a,b}$.
On the other hand, repeated action of $\SED$ generates the module $\mW_{a	,b;j}$ given by
$$
\mW_{a,b;j} = {\rm span}_{\mC} \{  \SE^{\dagger s} H^S_{a,b;j} : s = 0,1,2, \ldots, a-b   \}
$$
Each of the spaces $\mW_{a,b;j}$ is a realization of a finite dimensional irreducible $\gsl(2,\mC)$--module, which we denote by  $\mI_{a,b}$.
Finally, the space of $(a,b)$--homogeneous symplectic harmonic polynomials $\HS_{a,b}$ is a realization of the irreducible Sp$( p)$--module with highest weight $(a,b)_s = (a+b,0,\ldots,0)_s$ of length $p$, which we denote by $\mH_{a,b}$.
For all $(a,b)$ with $a>b$ the tensor product
$$
\left( \mI^\infty_{a,b} \otimes \mI_{a,b} \right) \otimes \mH_{a,b}
$$
then is an irreducible $(\gsl(2,\mC) \oplus  \gsl(2,\mC) ) \times \gsp_{2p}(\mC)$--module. When regarded as a 
$\gso(4,\mC)$--module it contains as many copies of $\mI^\infty_{a,b} \otimes \mI_{a,b}$ as the dimension of $\mH_{a,b}$; when regarded as an $\gsp_{2p}(\mC)$--module it contains infinitely many copies of $\mH_{a,b}$. The symplectic Fischer decompositions (\ref{symplectic-decomp})(\ref{symplectic-decomp-bis})(\ref{symplectic-decomp-ter}) may thus be reformulated as follows.

\begin{theorem}
Under the joint action of $(\gsl(2,\mC) \oplus  \gsl(2,\mC) ) \times {\rm Sp}( p)$, the space of complex valued polynomials $\mcP(\mR^{4p};\mC)$ is isomorphic with the multiplicity free irreducible direct sum decomposition
$$
\bigoplus_{a \geq b=0}^\infty \ \left( \mI^\infty_{a,b} \otimes \mI_{a,b} \right) \otimes \mH_{a,b}
$$
where $\mI^\infty_{a,b}$ is a Verma $\gsl(2,\mC)$--module with lowest weight $a+b+2p$, $\mI_{a,b}$ is an irreducible $\gsl(2,\mC)$--module with highest weight $a-b$ and $\mH_{a,b}$ is an irreducible $\gsp_{2p}(\mC)$--module with highest weight $(a+b,0,\ldots,0)$.
\end{theorem}


\end{document}